\documentclass[11pt]{article}
\usepackage{amsmath, amssymb, amsfonts}
\usepackage{color, soul}
\usepackage{epsfig}

\def\be{\begin{equation}}
\def\ee{\end{equation}}
\def\bea{\begin{eqnarray}}
\def\eea{\end{eqnarray}}
\def\bes{\begin{eqnarray*}}
\def\ees{\end{eqnarray*}}

\def\<{\langle}
\def\>{\rangle}

\title{The minimal periodic solutions for superquadratic autonomous Hamiltonian systems without the Palais-Smale condition}
\author{
Yuming Xiao \thanks{Supported by the NNSF, No. 11871356,
 e-mail: yumingxiao@scu.edu.cn.}\\
 School of Mathematics, Sichuan University, \\Chengdu 610064,   China \\
 Gaosheng Zhu \thanks{Supported by the NNSF, No. 11871368, e-mail: gaozsc@163.com, Corresponding author. }\\
 School of Mathematics, Harbin Institute of Technology, \\Harbin 150001,   China \\
}

\date{}
\begin{document}
\newtheorem{definition}{Definition}[section]
\newtheorem{theorem}{Theorem}[section]
\newtheorem{lemma}{Lemma}[section]
\newtheorem{corollary}{Corollary}[section]
\newtheorem{example}{Example}[section]
\newtheorem{property}{Property}[section]
\newtheorem{proposition}{Proposition}[section]
\newtheorem{remark}{Remark}[section]

\newcommand{\qed}{\nolinebreak\hfill\rule{2mm}{2mm}
\par\medbreak}
\newcommand{\Proof}{\par\medbreak\it Proof: \rm}
\newcommand{\rem}{\par\medbreak\it Remark: \rm}
\newcommand{\defi}{\par\medbreak\it Definition : \rm}
\renewcommand{\thefootnote}{\arabic{footnote}}

\maketitle
\begin{abstract}
{\it In this paper, we prove the existence of periodic solutions with any prescribed minimal period $T>0$
 for even second order Hamiltonian systems and convex first order Hamiltonian systems under the weak Nehari
  condition instead of Ambrosetti-Rabinowitz's. To this end, we shall develop the method of Nehari manifold to directly deal with a
frequently occurring problem where the Nehari manifold is not a manifold.

   }
\end{abstract}

{\bf Key words}:   Hamiltonian systems; Periodic solutions; Minimal period; Nehari manifold.

{\bf AMS Subject Classification}: 53C22, 58E05, 58E10.

\renewcommand{\theequation}{\thesection.\arabic{equation}}
\makeatletter\@addtoreset{equation}{section}\makeatother

\section{Introduction and the main results}

 Consider the following first order autonomous Hamiltonian system
\begin{equation}
\label{1stHS}
\dot{z}=JH'(z),\quad \forall z\in\mathbb{R}^{2n},
\end{equation}
where $n$ is a positive integer,
$J=\bigg( \begin{matrix}
    0 & I \\
    -I & 0
 \end{matrix} \bigg)$ is the standard symplectic matrix with $I$ the $n\times n$ identity matrix,
 $H:\mathbb{R}^{2n}\to\mathbb{R}$ is a function of class $C^1$ and $H'$ denotes the gradient of $H$.

In his seminal work \cite{Rab1978} of 1978, Rabinowitz introduced the variational formulation for
the system (\ref{1stHS}),
$$
f(z)={1\over2} \int_{0}^{T}(J\dot{z},z)dt-\int_{0}^{T}H(z)dt,\quad\forall z\in H^{1\over2}(S_{T},\mathbb{R}^{2n}),
$$
where $S_{T}=\mathbb{R}/(T\mathbb{Z})$ (cf. \S 3.1 of Abbondandolo\ \cite{Abb2001}).
He then used the variational method to establish the existence of non-constant periodic solutions with any
prescribed period $T>0$, if $H\geq0$ is superquadratic at the origin and satisfies
the well-known Ambrosetti-Rabinowitz condition:
\begin{equation}
\label{AR}
\exists \mu>2,\ \bar{r}>0,\ s.t.,\ \forall |z|\geq \bar{r},\ 0<\mu H(z)\leq H'(z)z.
\end{equation}
Because a $T/k$-periodic solution of (\ref{1stHS}) with $k\notin \mathbb{N}\backslash\{1\}$  can be also seen as a $T$-periodic one,
Rabinowitz further conjectured that (\ref{1stHS}) possesses a non-constant solution with any prescribed minimal period $T$.

This conjecture is still open now. However, there have been a great deal of works on it during the past over forty years.
In 1981, Ambrosetti and Mancini  \cite{AM1981} used the Clarke dual action principle and the method of Nehari manifold to
get solutions of minimal period for a class of strictly convex Hamiltonian systems.  Their idea was generalized  by Girardi and Matzeu \cite{GM1983}
of 1983 to deal with the non-convex Hamiltonian systems, by adding some other conditions to ensure that 0 is an isolated critical value of the direct action functional
(cf. also  Girardi and Matzeu \cite{GM1986} of 1986).
In 1985, a significant progress was made by Ekeland and Hofer \cite{EH1985} where they used the Clarke dual action principle, Ekeland's index theory and
 Hofer's topological characterization of the mountain pass point to prove Rabinowitz's conjecture for the strictly convex Hamiltoian systems.
 Later in 1990, Coti-Zelati, Ekeland and Lions \cite{CEL1990} proved the similar result as that of \cite{EH1985}, with (\ref{AR}) replaced by
 $ \lim_{|z|\to+\infty}H''(z)=+\infty. $  We refer the reader to Chapter IV of Ekeland \cite{Eke1990} for discussions and references before 1990.

  In 1997, Dong and Long \cite{DL1997} developed the Maslov index iteration formulae of symplectic paths and
  gave an intrinsic explanation and an extension to Ekeland-Hofer's theorom,
  where the positive assumption of $H''(z)$ in \cite{EH1985} is replaced by the semi-positive one together with
  another type of positive condition (cf. also Liu and Long \cite{Liu-Long2000} of 2000).
  Based on their work, Fei, Kim and Wang \cite{FKW1999} of 1999 obtained a new relationship between the iteration number and the Maslov indices and
  proved under the semi-positive assumption that (\ref{1stHS}) has at least a $T$  periodic solution  with $T/k$ its minimal period for some $k$ satisfying $1\leq k\leq 2n.$
  We refer the reader to Chapter 13 of Long \cite{Long2002} along this line before 2002.
 Notice also that in addition to the semi-positive assumption, if $H$ is symmetric (i.e., $H(-z)=H(z)$),
 Zhang \cite{Zhang2014} of 2014 investigated the iteration formulae of symmetric Maslov indices
  and gave a confirming answer to Rabinowitz's conjecture.

  There is also a version of Rabinowitz's conjecture for the second order Hamiltonian system
  \begin{equation}
  \label{2ndHS}
  \ddot{x}+V'(x)=0,\quad\forall x\in\mathbb{R}^{N}.
  \end{equation}
  Indeed, many works mentioned above considered the second order case too and obtained similar results as those for the first order one, such as \cite{CEL1990}  and \cite{Zhang2014}.
  Nevertheless, the first result under precise Rabinowitz's superquadratic conditions is due to Long \cite{Long1994} of 1994,
  where he established the iteration formulae of the symmetric Morse index
  and proved by the saddle point theorem that (\ref{2ndHS}) possesses at least a $T$-periodic even solution with minimal period $T/k$ for some integer $k$ satisfying $1 \leq k\leq N+2$
  which was subsequently improved by himself to  $1\leq k \leq N+1$  in \cite{Long1997} of 1997. Notice that it was also proved $k=1$ or $2$  in \cite{Long1994}, if moreover $V''(x)$ is positive for $x\in\mathbb{R}^{N}\backslash\{0\}$.

     As for the second order even Hamiltonian systems (i.e., $V(x)=V(-x)$), by reducing the working space to a smaller one so that the mountain pass theorem can be applied,
     Long \cite{Long1993} of 1993  proved a similar result with $k=1$ or $3$ as that of \cite{Long1994}. If moreover $V''(x)$ is semi-positive, Fei, Kim and Wang \cite{FKW2001} of 2001
     eliminated the case of $k=3$ and Zhang \cite{Zhang2014} gave a different proof of such a result. Quite different from the method of the index iteration,
     in 2010 Xiao \cite{Xiao2010} used the method of Nehari manifold to get a $T$-minimal periodic solution for such a problem,
      under the following condition
 \begin{equation}
\label{Xiao1}
\exists\theta>1,\ s.t.,\ 0<\theta V'(x)\cdot x\leq V''(x)x\cdot x,\quad\forall x\in\mathbb{R}^{2n}\backslash\{0\},
\end{equation} which was later weakened to some extent by Xiao and Shen \cite{Xiao-Shen2020} of 2020.

   For more works and recent progress on Rabinowitz's conjecture for both (\ref{1stHS}) and (\ref{2ndHS}),
    we refer the reader to the excellent review paper by Long \cite{Long2021} of 2021.

In this paper, we first prove a result on the existence of minimal period solution for the  second order Hamiltonian system
(\ref{2ndHS}).
\begin{theorem}
\label{mainresult1}
Suppose that $V\in C^{2}(\mathbb{R}^{N},\mathbb{R})$ satisfies the following conditions

${\rm (V1)}$ $\lim_{x\to0}{V(x)\over|x|^{2}}=0$,

${\rm (V2)}$ $\lim_{x\to\infty}{V(x)\over|x|^{2}}=+\infty,$

${\rm (V3)}$   $0< V'(x)\cdot x \leq V''(x)x\cdot x,\ \forall x\in\mathbb{R}^{N}\backslash\{0\},$

${\rm (V4)}$  $V(x)=V(-x)$.

\noindent Then for every $T>0$,  there is a $T$-periodic solution of (\ref{2ndHS}) with $T$ its minimal period,
which is even about $t = 0,\ T/2,$ and odd about $t = T/4,\ 3T/4.$
\end{theorem}
\begin{remark}
 (V3) is not comparable with the Ambrosetti-Rabinowitz condtion (\ref{AR}). In some sense it can be
considered weaker because it allows $V$ to grow quadratically, although in our context $V$ should  grow  eventually  more than quadratic due to (V2).
A typical class of such examples is  $V(x)=|x|^{2}(\ln(1+|x|^{p}))^{q}$  with $p,q>0$, which satisfy (V1)-(V4), but not (\ref{AR}).
\end{remark}
\begin{remark}
Conditions (V1), (V2) and (V3) can be derived from (\ref{Xiao1}),
provided that $V(0)=0$ which is not essential. Indeed, (\ref{Xiao1}) is also stronger than (\ref{AR}), cf. \cite{Xiao2010}.
\end{remark}

Our next result is on the existence of minimal period solution for the  first order Hamiltonian system (\ref{1stHS}).

\begin{theorem}
\label{mainresult2}
Suppose  $H\in C^{1}(\mathbb{R}^{2n},\mathbb{R})$ is strictly convex and satisfies the following conditions.

{\rm(H1)} $\lim_{x\to0}{H(x)\over|x|^{2}}=0$.

{\rm(H2)}  There are $\beta>2$  and positive constants $r$, $a_{1}$ and $a_{2}$ such that
\begin{equation}
\label{crucialcondition1}
a_{1}|x|^{\beta}\leq H(x)\leq a_{2}|x|^{\beta},\quad\forall|x|\geq r.
\end{equation}

{\rm(H3)} For every $y\in\mathbb{R}^{2n}\backslash\{0\}$, $s\mapsto{G'(sy)y\over s}$ is non-increasing about $s>0,$
where $G=H^{*}$ is the Fenchel transform of $H$.

Then, for every $T>0$, there is a $T$-periodic solution of (\ref{1stHS}) with $T$ its minimal period.
\end{theorem}
\begin{remark}As far as we know, very few works on the minimal period problem discuss the case when $H$ is only of class $C^1$.
Nevertheless, if $H\in C^{2}(\mathbb{R}^{2n},\mathbb{R})$ and $H''(x)$ is positive definite for $x\neq0$,
then $G\in C^{2}(\mathbb{R}^{2n}\backslash\{0\},\mathbb{R})$ and {\rm(H3)} is equivalent to the (dual) weak Nehari condition
$$G''(y)y\cdot y\leq G'(y)y,\quad\forall y\in\mathbb{R}^{2n}\backslash\{0\}.$$
\end{remark}

This paper is organized as follows. In section 2, we develop the method of Nehari manifold  to directly deal with a
frequently occurring problem where the Nehari manifold is not a manifold, and prove Theorem \ref{abstract1} and Theorem \ref{abstract2} which are independently interesting.
 Then  in section 3, we apply Theorem \ref{abstract1} to prove  Theorem \ref{mainresult1}.
Finally in section 4, by combining Theorem \ref{abstract2} with the Clarke dual action principle we give a proof of Theorem \ref{mainresult2}.
\section{A generalization of the method of Nehari manifold}
The method of Nehari manifold  is widely used to find ground state solutions of many partial differential equations (cf. Szulkin and Weth \cite{Szulkin-Weth2010} of 2010).
Let $E$ be a real Banach space and $\Phi\in C^{1}(E,\mathbb{R})$. Define
$$\mathcal{N}=\{x\in E\backslash\{0\}\mid \Phi'(x)x=0\}.$$
It is obvious that if $x\neq0$ is a critical point of $\Phi$, it must lie on $\mathcal{N}$.
Thus we need only devote ourselves to finding critical points of $\Phi$ on  $\mathcal{N}$
and this is the basic idea of the method of Nehari manifold.

In the literature, $\mathcal{N}$ is called the Nehari manifold of $\Phi$,
although it may not be a manifold in general.  Nevertheless,  to make this method work effectively,
some strict conditions are required to ensure $\mathcal{N}$ to be a manifold of at least class $C^0$ in most applications
 (cf. \cite{Liu-Wang2004}, \cite{Szulkin-Weth2010}, and \cite{Xiao2010}).

If  $\mathcal{N}$ is of class $C^1$ fortunately, then we can carry out the critical point theory to $\Phi$ on $\mathcal{N}$.
One of the advantage is that $\Phi\big|_{\mathcal{N}}$ may satisfy the Palais-Smale condition under weaker conditions than $\Phi\big|_{E}$.
However, if $\mathcal{N}$ is only of class $C^0$ or not a manifold at all, the above way  fails since $(\Phi\big|_{\mathcal{N}})'$ is not well defined.
To overcome this difficulty, some indirect arguments are needed (cf. \cite{PKS2017}, \cite{Szulkin-Weth2010}, \cite{Tang2015} and \cite{Xiao-Shen2020}).

This paper seems to be the first one that directly deals with a frequently occurring problem where the Nehari manifold is not a manifold.
The key ingredients of our arguments are completely different from and definitely much simpler and more direct than the above mentioned works.
In fact, our method is more topological rather than analytic as usual. Thus it does not require any type of compact conditions, such as the Palais-Smale condition.

The first result of this section can be  formulated as follows.
 \begin{theorem}
 \label{abstract1}
 Let $E$ be a real reflexive Banach space with $||\cdot||$  its  norm, and $$\Phi(x)={1\over2}||x||^{2}-I(x)$$
 where  $I\in C^{1}(E,\mathbb{R})$ satisfies the following conditions

 {\rm($i$)} $I(x)$ is weakly continuous;

 {\rm($ii$)} $\lim_{x\to0}{I(x)\over ||x||^{2}}=0;$

 {\rm($iii$)} ${I(sx)\over s^{2}}\to+\infty$, uniformly for $x$ on weakly compact subsets of $E\backslash\{0\}$ as $s\to+\infty$;

 {\rm($iv$)} $s\mapsto {I'(sx)x\over s}$ is non-decreasing for all $x\neq0$ and $s>0$;

 Then, there is a critical point $\bar{x}$ of $\Phi$ such that
 \begin{equation}
 \label{0809a}
 \Phi(\bar{x})=\max_{s\geq0} \Phi(s\bar{x})=\inf_{||e||=1}\max_{s\geq0} \Phi(se)=\inf_{x\in\mathcal{N}}\Phi(x)>0.
 \end{equation}
 \end{theorem}

 Theorem \ref{abstract1} considerably improves Theorem 12 of \cite{Szulkin-Weth2010}, if only the
conclusion on the existence is concerned. The novelty of our proof is that we use not only the Nehari manifold but also the level set
to detect the critical points characterized by (\ref{0809a}).  From the view of result, it can be regarded as an application of Morse theory
for $\Phi|_{E}$ without the Palais-Smale condition. We also mention that our idea can be  used to study the Nehari-Pankov
manifold, namely the generalized Nehari manifold (cf. \cite{Szulkin-Weth2009} and Chapter 4 of \cite{Szulkin-Weth2010}).
We will give  more comments on it and explain how to use it to study the minimal period problem for general Hamiltonian systems
without the even or convex  assumption at the end of section 4.

 Our next result, in some sense, can be regarded as a dual version of Theorem \ref{abstract1}. To state it, we first provide some notations.
 Let $E$ be a real Banach space  with $||\cdot||$ its norm  and $a:E\times E\to\mathbb{R}$ is a continuous, symmetric, bilinear form.
 Then, the set
 $$\mathcal{P}^{0}=\{u\in E\mid a(u,u)=0\}$$
 divides $E$ into two parts, $\mathcal{P}^{+}=\{u\in E\mid a(u,u)>0\}$ and $\mathcal{P}^{-}=\{u\in E\mid a(u,u)<0\}$.
 Let $\mathcal{S}=\{u\in E\mid||u||=1\}$ be the unit sphere of $E$ and
  denote by $\mathcal{S}^{+}=\mathcal{S}\cap\mathcal{P}^{+}$, $\mathcal{S}^{-}=\mathcal{S}\cap\mathcal{P}^{-}$ and $\mathcal{S}^{0}=\mathcal{S}\cap\mathcal{P}^{0}$.
 \begin{theorem}
\label{abstract2}
Let $E$ be a real reflexive Banach space with $||\cdot||$ its norm  and $\Phi$ is a $C^{1}$-functional on $E$ with the form
 $$\Phi(u)={1\over2}a(u,u)+b(u),$$
where $a:E\times E\to\mathbb{R}$ is a continuous, symmetric, bilinear form with nonempty negative eigenspace. Suppose that

{  (i) }$a(u,u)$ is weakly continuous and $b(u)$ is weakly lower semi continuous;

{ (ii)} $b(0)=0$, $\lim_{u\to0}{b(u)\over||u||^{2}}=+\infty $ and $\lim_{u\to\infty}{b(u)\over||u||^{2}}=0;$

{ (iii)} $b'(u)u\geq b(u)>0$ for $u\neq0$;

{ (iv)} for every $u\neq 0$, $s\mapsto{b'(su)u\over s}$ is  non-increasing on $(0,+\infty)$.

Then, $\Phi$ has a critical point $\bar{u}\in\mathcal{P}^{-}$ with
\begin{equation}
\label{0809b}
\Phi(\bar{u})=\max_{s\geq0}\Phi(s\bar{u})=\inf_{e\in \mathcal{S}^{-}}\max_{s\geq0}\Phi(se)=\inf_{u\in\mathcal{N}}\Phi(u)>0.
\end{equation}
\end{theorem}

\begin{remark}
 A trivial $b$ which satisfies all the conditions of Theorem \ref{abstract2} is $||u||^{\alpha}$ with $1<\alpha<2.$
 For  application to concrete problems (cf. Ekeland \cite{Eke1990}, p.108), if necessary, many of the conditions in Theorem \ref{abstract2} can be relaxed,
  but they are enough for our need.
\end{remark}
\begin{remark} Generally, the topologies of the Nehari manifolds in Theorem \ref{abstract1} and Theorem \ref{abstract2} are different. Roughly speaking,
the former is a thickened sphere, while the latter is a thickened revolving paraboloid.
\end{remark}

In this section, we prove Theorem  \ref{abstract1} in details. The proof of Theorem \ref{abstract2} is similar,
so we only sketch it and stress the differences between them. First, we investigate the topology of the Nehari manifold in Theorem  \ref{abstract1}
$$\mathcal{N}=\{x\in E\backslash\{0\}\mid\Phi'(x)x=0\},$$
and the result can be stated as follows.
\begin{lemma}
 \label{topology}
 Under the assumptions of Theorem \ref{abstract1}, the Nehari manifold
 \begin{equation}
 \label{20230217A}
 \mathcal{N}=\bigcup_{||e||=1}\mathcal{N}_{e},
 \end{equation}
 where $\mathcal{N}_{e}$ is a point or a bounded line segment which lies in the ray $\overrightarrow{oe}$ and is isolated from the origin.
  Moreover, $\mathcal{N}$ divides $E$ into two path-connected components.
 \end{lemma}
 \Proof By ($ii$), there exists sufficiently small $\rho_{0}>0$ such that
 \begin{equation}
 \label{20230304a}
 \Phi(x)\geq{1\over 2}||x||^{2}-{1\over4}||x||^{2}={1\over4}\rho^{2},\quad\forall||x||=\rho\leq\rho_{0}.
 \end{equation}
 Now fix arbitrarily $e\in E$ with $||e||=1$ and define
 $$\phi(s)=\Phi(se)={1\over2}s^{2}-I(se),\quad\forall s\geq0.$$
 According to ($iii$), it holds
\begin{equation}
\label{20230217d}
\lim_{s\to+\infty}\phi(s)=\lim_{s\to+\infty}s^{2}\left({1\over2}-{I(se)\over s^{2}}\right)=-\infty.
\end{equation}
Then $\phi(0)=0$, together with (\ref{20230304a}) and (\ref{20230217d}), yield the existence of $\bar{s}>0$ satisfying
 \begin{equation}
 \label{20230304b}\phi(\bar{s})=\max_{s\geq0}\phi(s)>0,\quad\phi'(\bar{s})=0.
 \end{equation}

Notice that
$$\begin{aligned}
\phi'(s)&=s-I'(se)e =s\left(1-{I'(se)e\over s} \right),\quad\forall s>0,
\end{aligned}$$ and  so the zeroes of $\phi'(s)$ is the same as those of $1-{I'(se)e\over s}$.
 Due to $(iv)$,  the latter function is non-increasing  about $s>0 $.
Therefore, the zeroes of $\phi'(s)$ $(s>0)$ is at most a bounded and connected interval $I$.
It then follows from (\ref{20230304b}) that $\bar{s}\in I\subset(0,+\infty)$ and
$$\phi'(s)=0,\quad \phi(s)=\phi(\bar{s})=\max_{s\geq0}\phi(s),\quad\forall s\in I.$$

 Again by (\ref{20230304a}), we have $$c_{E}=\inf_{||e||=1}\max_{s\geq0}\Phi(se)>0.$$
  Since $\lim_{x\to 0}\Phi(x)=0$, $\mathcal{N}$ is isolated from the origin and the previous conclusion is proved.

Now we come to prove the latter conclusion. Define
$$\Omega_{1}=\{x\in E\backslash\mathcal{N}\mid\overline{ox}\cap\mathcal{N}=\emptyset\},\quad\Omega_{2}=\{x\in E\backslash\mathcal{N}\mid\overline{ox}\cap\mathcal{N}\neq\emptyset\},$$
where $\overline{ox}$ is the line segment connecting the origin $o$ and $x$.

It is obvious that $\Omega_{1}$ is path-connected, since every point $x\in\Omega_{1}$ can be connected with the origin by $\overline{ox}\subset\Omega_{1}$.
Now let $x_{1}$ and $x_{2}$ be two points in $\Omega_{2}$. Choose a plane $\Pi$ containing the lines $ox_{1}$ and $ox_{2}$.
Observing that $\Pi\cap \mathcal{N}$ is a bounded set, we define
$K=\max\{||x||,x\in\Pi\cap \mathcal{N}\}$ and  choose
 $$R>\max\{||x_{1}||,||x_{2}||,K\}.$$
 Then we can connect $x_{1}$ with $x_2$ by the path: first $\{sx_{1}\mid 1\leq s\leq {R\over||x_{1}||}\}$,
 then from $R {x_{1}\over||x_{1}||}$ to $R {x_{2}\over||x_{2}||}$ along the circle $\{x\in\Pi,\ ||x||=R\}$, finally $\{ sx_{2} \mid 1\leq s\leq {R\over||x_{2}||}\}^{-1}$.
  Thus $\Omega_{2}$ is also path-connected.

 To prove $\Omega_{1}$ and $\Omega_{2}$ are separated, consider the functional
$$\Phi'(x)x=||x||^{2}-I'(x)x=||x||\cdot\phi'(s)\bigg|_{s=||x||}.$$
Thus, $\Phi'(x)x>0$ on $\Omega_{1}\backslash\{0\}$, $\Phi'(x)x=0$ on $\mathcal{N}$, and $\Phi'(x)x<0$ in $\Omega_2$.
Since $\Phi'(x)x$ is continuous, $\Omega_{1}$ and $\Omega_{2}$ are  separated by $\mathcal{N}$.\hfill$\Box$
\begin{remark}
\label{intuitive}
Intuitively, along every ray from the origin, $\Phi$ is strictly increasing from $0$ until it attains its maximum value.
Then it may remain constant in a bounded interval. After that, it shall strictly decrease to $-\infty.$
\end{remark}
\begin{remark}
\label{topexample}
In general, $\mathcal{N}$ is not a manifold. One may suspect that $\mathcal{N}$ contains a manifold
 $\mathcal{N}_{1}$ of codimension 1 which divides $E$ into two path-connected components.
However, we do not know whether such a phenomenon appears as the following simple example.

 Let
$$\begin{aligned}
\Omega=&[-2,2]\times[-2,2]\\
\mathcal{N}=&\left\{(x,y)\in\Omega\mid y=\sin\Big({1\over x}\Big),\ x\in[-2,2]\backslash\{0\}\right\}\cup\{(0,y)\mid -1\leq y\leq 1\},\\
 \Omega_{1}=&\left\{(x,y)\in\Omega\mid y<\sin\Big({1\over x}\Big)\right\},\quad\Omega_{2}=\left\{(x,y)\in \Omega\mid y>\sin\Big({1\over x}\Big)\right\}.
 \end{aligned}
 $$
It is obvious that $\mathcal{N}$ divides $\Omega$ into two path-connected components,
but there is no such a manifold $\mathcal{N}_{1}\subset\mathcal{N}$ which possesses the same property.
\end{remark}

Now we come to give a proof of Theorem \ref{abstract1}.

{\bf Proof of Theorem \ref{abstract1}: }

Define $$c_{E}=\inf_{||e||=1}\max_{s\geq0}\Phi(se).$$
Due to (\ref{20230304a}), it holds that $c_{E}>0.$ Take arbitrarily a sequence $\{x_{k}\}\subset\mathcal{N}$ such that
$$\lim_{k\to+\infty}\Phi(x_{k})= c_{E}.$$

{\bf Claim 1:} There is a subsequence of $\{x_{k}\}$ weakly converging to some $\bar{x}\in E\backslash\{0\}$.

 We first prove that $\{x_{k}\}$ is bounded. Passing to a subsequence if necessary, we may assume  $||x_{k}||\to+\infty$ by contradiction.
Let  $y_{k}={x_{k}\over||x_{k}||}$ and suppose $y_{k}\rightharpoonup \bar{y}$ without loss of generality.
If $\bar{y}\neq0$, then by ($iii$) we have
$$0=\lim_{k\to+\infty}{\Phi(x_{k})\over ||x_{k}||^{2}}={1\over 2}-\lim_{k\to+\infty}{I(x_{k})\over||x_{k}||^{2}}={1\over 2}-\lim_{k\to+\infty}{I(||x_{k}||y_{k})\over||x_{k}||^{2}}=-\infty,$$
which is obviously absurd. While if $\bar{y}=0$, then for  arbitrary $R>0$ it follows from ($i$) that
\begin{equation}
\label{20230325a}
c_{E}=\lim_{k\to+\infty}\Phi(x_{k})\geq\lim_{k\to+\infty}\Phi(Ry_{k})={1\over2}R^{2} -\lim_{k\to+\infty}I(Ry_{k})={1\over2}R^{2}.
\end{equation}
Thus, if we choose $R>\sqrt{2 c_{E} }$ a priori, a contradiction is obtained.

Assume without loss of generality that $\{x_{k}\}$ weakly converges to $\bar{x}$.
If $\bar{x}=0$, then by the geometric characteristic of $x_{k}$ and ($i$)
$$c_{E}=\lim_{k\to+\infty}\Phi(x_{k})\geq\lim_{k\to+\infty}\Phi(Rx_{k})
={1\over2}R^{2}\lim_{k\to+\infty}||x_{k}||^{2}\geq{1\over2}d^{2}R^{2},$$
where $d=\inf_{\mathcal{N}}||x||>0$ due to Lemma \ref{topology}. Thus, if we choose $R> {1\over d}\sqrt{2 c_{E}}$, a contradiction is obtained and Claim 1 is proved.

{\bf Claim 2:} $\max_{s\geq0}\Phi(s\bar{x})=c_{E}.$

Indeed, by the topology of $\mathcal{N}$ and the definition of $c_E$, there is $\alpha>0$ such that
\begin{equation}
\label{20230228a}
\Phi(\alpha\bar{x})=\max_{s\geq0}\Phi(s\bar{x})\geq c_{E}.
\end{equation}

On the other hand, since ${1\over2}||x||^{2}$ is convex and continuous, it is weakly lower semi continuous by Theorem 1.2 of \cite{MW1989}.
Then it follows from ($i$) that $\Phi(x)={1\over2}||x||^{2}-I(x)$ is a weakly lower semi continuous functional.
As a result, we get
\begin{equation}
\label{20230228b}
\Phi(\alpha\bar{x})\leq{\lim\inf}_{k\to+\infty}\Phi(\alpha x_{k})\leq {\lim\inf}_{k\to+\infty}\Phi(x_{k})= c_{E}.
\end{equation}
The inequalities (\ref{20230228a}) and (\ref{20230228b}) then yield Claim 2.

{\bf Claim 3:} There is a critical point $z$ of $\Phi\big|_{E}$  such that
\begin{equation}
\label{20230326A}
\Phi(z)=\max_{s\geq0}\Phi(sz)=c_{E}.
\end{equation}

To this end, consider the level set
$$L(c_{E})=\{x\in E\mid \Phi(x)=c_{E}\}.$$
Define
$$\begin{aligned}
\Sigma_{+}=&\{e\in E\mid ||e||=1,\ \max_{s\geq0}\Phi(se)>c_{E}\},\\
\Sigma_{0}=&\{e\in E\mid ||e||=1,\ \max_{s\geq0}\Phi(se)=c_{E}\},\\
\mathcal{N}_{0}=&\bigcup_{e\in\Sigma_{0}}\mathcal{N}_{e}.\\
\end{aligned}$$
Notice that $\mathcal{N}_{0}\neq\emptyset$ by Claim 2 and it is bounded according to the proof of Claim 1.

For every $e\in\Sigma_{+}$, by Remark \ref{intuitive} there are exactly $x_{1}=x_{1}(e)\in\Omega_{1}$ and $x_{2}=x_{2}(e)\in\Omega_{2}$ such that
$$\Phi(x_{1})=\Phi(x_{2})=c_{E},\quad \Phi'(x_{1})x_{1}>0,\quad  \Phi'(x_{2})x_{2}<0.$$
By the implicit function theorem, $\mathcal{M}_{1}=\bigcup_{e\in\Sigma_{+}}x_{1}(e)$ and $\mathcal{M}_{2}=\bigcup_{e\in\Sigma_{+}}x_{2}(e)$ are
manifolds of codimension $1$ in $\Omega_{1}$ and $\Omega_{2}$ respectively, and
$$h:\mathcal{M}_{1}\to \mathcal{M}_{2},\quad x_{1}(e)\to x_{2}(e),$$
is a homeomorphism.
Furthermore, $\mathcal{M}_{1}\cup\mathcal{N}_{0}$ (resp. $\mathcal{M}_{2}\cup\mathcal{N}_{0}$)
divides $E$ into two path-connected components.
Assume by contradiction that Claim 3 is false. Then by the implicit function theorem, the level set
 $L(c_{E})=\mathcal{M}_{1}\cup\mathcal{N}_{0}\cup\mathcal{M}_{2}$
is a manifold of codimension $1$.
However, it can be checked directly from all the above topological information that $L(c_{E})$ can not be a manifold.

For this purpose, we denote by $\partial\mathcal{M}_{1}$ the boundary of ${\mathcal{M}}_{1}$.
Due to $\partial\mathcal{M}_{1}\subseteq\mathcal{N}_{0},$  the closure $\overline{\mathcal{M}}_{1}=\mathcal{M}_{1}\cup\partial\mathcal{M}_{1}$
is a submanifold of $L(c_{E})$ with codimension $1$, and  $\partial \overline{\mathcal{M}}_{1} $  is a union of some manifolds
 which are homeomorphic to spheres of codimension $2$.
Notice that $\partial \overline{\mathcal{M}}_{1}$ is a subset of $\partial {\mathcal{M}}_{1}$ in general.

If $A=\partial{\mathcal{M}}_{1}\backslash\partial \overline{\mathcal{M}}_{1}\neq\emptyset,$  then for every  $p\in A$,
there is an open neighborhood $U$ of $p$ with codimension $1$ such that  $U\subseteq \overline{\mathcal{M}}_{1}.$
Define $$\mathcal{N}(A)=\{\mathcal{N}_{a}\mid a\in A\}.$$
Thus if $L(c_{E})$ is a manifold of codimension $1$, then $\mathcal{N}(A)=A$ and so $A\subset\partial\mathcal{M}_{2}$.
It then yields that the intersection of a neighborhood of $A$ in $E$ with $\mathcal{M}_{2}$ is nonempty
which is obviously absurd.

Since the above arguments are valid if we replace $\mathcal{M}_{1}$ by $\mathcal{M}_{2},$
we now suppose that  $\partial{\mathcal{M}}_{1}=\partial \overline{\mathcal{M}}_{1}$ and $\partial{\mathcal{M}}_{2}=\partial \overline{\mathcal{M}}_{2}$.
 Choose arbitrarily a connected component  $A_{1}$ in $\partial{\mathcal{M}}_{1}$, which is a manifold homeomorphic to a sphere of codimension $2$.
 According to the above arguments, $A_{1}$ shall be connected by  $\mathcal{N}({A_1})$ to some manifold $A_{2}$ in $\partial\mathcal{M}_{2}$ with codimension $2$.
 Again by the fact that $L(c_{E})$ is a manifold of codimension $1$, the cylinder formed by $A_{1}$, $A_{2}$ and $\mathcal{N}(A_{1})$ is a manifold of codimension $1$.

 However, as we mentioned previously, $\mathcal{M}_{1}\cup\mathcal{N}_{0}$
 (resp. $\mathcal{M}_{2}\cup\mathcal{N}_{0}$) divides $E$ into two path-connected components.
 To close this cylinder, there has to be still a simply-connected subset of $\mathcal{N}_{0}$ which adjoins $A_{1}$ or $A_{2}$.
As a result, $L(c_{E})$ can not be a manifold because of some points in $A_{1}$ or $A_{2} $. This leads to a contradiction.
 \hfill$\Box$

\begin{remark}
\label{nonmfld}
 Claim 2 is crucial in the proof of Theorem \ref{abstract1}. Indeed, if it did not hold,
then $\mathcal{M}_{1}$ and  $\mathcal{M}_{2}$ would be homeomorphic to $S^{\infty}(a)$ and  $S^{\infty}(b)$ respectively,
with $0<a<b$. Thus from a geometric perspective, $\mathcal{N}_{0}\neq\emptyset$ plays the role in gluing $\mathcal{M}_{1}$ with $\mathcal{M}_{2}$,
so that  $L(c_{E})=\mathcal{M}_{1}\cup\mathcal{N}_{0}\cup\mathcal{M}_{2}$ can not be a manifold.
\end{remark}

It is natural to ask whether one can use Morse theory to prove Claim 3. To answer this question, notice first that for $\epsilon>0$ sufficiently small,
the sub-level set
$$
\begin{aligned}
\Phi^{c_{E}-\epsilon}&=\{x\in E\mid\Phi(x)\leq c_{E}-\epsilon\}\\
&=\{x\in \Omega_{1}\mid\Phi(x)\leq c_{E}-\epsilon\}\cup \{x\in \Omega_{2}\mid\Phi(x)\leq c_{E}-\epsilon\},
\end{aligned}
$$
and by Remark \ref{intuitive} it consists of two connected components  inside $\{x\in \Omega_{1}\mid\Phi(x)= c_{E}-\epsilon\}$
and outside $\{x\in \Omega_{2}\mid\Phi(x)= c_{E}-\epsilon\}$ respectively.

On the other hand, $\Phi^{c_{E}+\epsilon}$ is connected by Claim 2.
Thus, Morse theory implies that $c_E$ is a critical value of $\Phi$, provided that $\Phi$ satisfies the Palais-Smale condition in $E$ (cf. \S 1.3 in \cite{Chang1993}).
Moreover, it follows from Remark \ref{intuitive} again that such critical points must lie in $\mathcal{N}_{0}$.
However, as we have mentioned in section 1, the assumption on the Palais-Smale condition for $\Phi$ is not suitable for our applications.

In the proof of Theorem 2.1 in \cite{Liu-Wang2004}, Liu and Wang use an argument of mathematical analysis to detect such a critical point of $\Phi$.
The geometric meaning therein can be described as follows. For simplicity of writing, we suppose that $\mathcal{N}_{0}=\{z\}$. If $z$ is not a critical point,
then there is $v\in E$ such that $\Phi'(z)v<0.$ As a result, for $\epsilon,\eta>0$ sufficiently small, it holds that
$$(1-\epsilon)z+\eta v\in\Omega_{1},\quad(1+\epsilon)z+\eta v\in\Omega_{2},$$
and the path $tz+\eta v$ $(1-\epsilon\leq t\leq 1+\epsilon)$ shall intersect with $\mathcal{N}$ at some point $\bar{t} z+\eta v$. Then,
\begin{equation}
\label{dis}
\Phi(\bar{t} z+\eta v)< \Phi(\bar{t} z)\leq\Phi(z)=c_{E},
\end{equation}
which contradicts to the definition of $c_E$. However, if $\mathcal{N}$ is not a manifold, such a way fails.
The reason is that the non-negligible distance between the ends of $\mathcal{N}_{z}$ may destroy the inequality (\ref{dis}).

To overcome this difficulty, we investigate  instead the topology of the level set $L(c_{E})$.  Needless to say, our
method  works for the above case which can indeed be regarded as the simplest one.
It is easy to see that in this case $L(c_{E})$ is homeomorphic to $S^{\infty}\vee S^{\infty}$ where $\vee$ is the wedge sum and so it is not a manifold.

Now we sketch a proof of Theorem \ref{abstract2}.

{\bf Proof of Theorem \ref{abstract2}:}

By $(ii)$, $\Phi(0)=0$ and $0$ is a local minimal point of $\Phi$.  Fix arbitrarily $e\in\mathcal{S}$ and consider
$$\phi(s)=\Phi(se)={1\over2}a(e,e)s^{2}+b(se),\quad\forall s\geq0.$$
Then,
$$\phi'(s)=a(e,e)s+b'(se)e=s\left(a(e,e)+{b'(se)e\over s}\right),\quad\forall s> 0.$$

If $e\in\mathcal{S}^{-}$,  it follows by $(ii)$ that $\lim_{s\to+\infty}\phi(s)=-\infty$. Therefore,
 there is $\bar{s}\in(0,+\infty)$ such that $\phi'(\bar{s})=0.$  By the same arguments as in Theorem \ref{abstract1}, the set of the critical points of $\phi$ on $(0,+\infty)$
 is a single point or a bounded and connected interval, which shares the same geometric meaning as we stated in Remark \ref{intuitive}.
 While if $e\in\mathcal{S}^{+}\cup \mathcal{S}^{0}$, then by $(iii)$
 $$\phi'(s)\geq b'(se)e>0,\quad\forall s>0.$$ As a result, the Nehari manifold of $\Phi$ is
  $\mathcal{N}=\bigcup_{e\in\mathcal{S}^{-}}\mathcal{N}_{e},$
 with $\mathcal{N}_{e}$ a single point or a bounded and connected interval. Again by $(ii)$, $\mathcal{N}$ is isolated from the origin.

 Define $$c_{E}=\inf_{e\in\mathcal{S}^{-}}\max_{s\geq0}\Phi(se)>0,$$
and take arbitrarily a sequence $\{u_{k}\}\subset\mathcal{N}$ such that
$$\lim_{k\to+\infty}\Phi(u_{k})= c_{E}.$$

 {\bf Claim 1:} There is a subsequence of $\{u_{k}\}$ weakly converging to some $\bar{u}\in E\backslash\{0\}$.

 Since $E$ is a reflexive Banach space, it is sufficient to prove that $\{u_{k}\}$ is bounded.
 If not, passing to a subsequence if necessary, we may assume $\lim_{k\to+\infty}||u_{k}||=+\infty.$
 Let $v_{k}={v_{k}\over||u_{k}||}$ and  suppose without loss of generality that $v_{k}\rightharpoonup v$. Then by ($i$) and ($ii$) we get
 $$0=\lim_{k\to+\infty}{\Phi(u_{k})\over||u_{k}||^{2}}=\lim_{k\to+\infty}{1\over2}a(v_{k},v_{k})+\lim_{k\to+\infty}{b(u_{k})\over||u_{k}||^{2}}={1\over2}a(v,v).$$
 Take $R>0$ arbitrarily, then
  $$c_{E}=\lim_{k\to+\infty}\Phi(u_{k})\geq\lim_{k\to+\infty}\Phi(Rv_{k})={1\over2}R^{2}\lim_{k\to+\infty}a(v_{k},v_{k})+\lim_{k\to+\infty}b(Rv_{k})=\lim_{k\to+\infty}b(Rv_{k}).$$
 Notice that ($ii$) and ($iii$) imply that there exist positive constants $\alpha$ and $\delta$ such that $||b(u)||\geq\alpha||u||$ for $||u||\geq\delta.$
 Thus if we choose $R>\max\{\delta,{c_{E}\over\alpha}\}$ a priori, a contradiction is obtained.

{\bf Claim 2:} $\max_{s\geq0}\Phi(s\bar{u})=c_{E}.$

Assume without loss of generality that $u_{k}\rightharpoonup\bar{u}$. Due to $(ii)$ and $(iii)$, there exists $\rho>0$ such that
 $a(u_{k},u_{k})=-b'(u_{k})u_{k}\leq-\rho<0.$    Then by $(i)$,
$$a(\bar{u},\bar{u})=\lim_{k\to+\infty}a(u_{k},u_{k})\leq-\rho<0.$$
According to  the topology of $\mathcal{N}$ and the definition of $c_E$, there is $\alpha>0$ such that
\begin{equation}
\label{20230228aa}
\Phi(\alpha\bar{u})=\max_{s\geq0}\Phi(s\bar{u})\geq c_{E}.
\end{equation}

On the other hand,  ($i$) implies that $\Phi(u)$ is a weakly lower semi-continuous functional and so
\begin{equation}
\label{20230228bb}
\Phi(\alpha\bar{u})\leq{\lim\inf}_{k\to+\infty}\Phi(\alpha u_{k})\leq {\lim\inf}_{k\to+\infty}\Phi(u_{k})= c_{E}.
\end{equation}
The inequalities (\ref{20230228aa}) and (\ref{20230228bb}) then yield Claim 2.

{\bf Claim 3:} There is a critical point $z$ of $\Phi\big|_{E}$ in $\mathcal{P}^{-}$  such that
\begin{equation}
\label{20230326AA}
\Phi(z)=\max_{s\geq0}\Phi(sz)=c_{E}.
\end{equation}

Similarly, we consider the level set
$$\begin{aligned}
L(c_{E})&=\{x\in E\mid \Phi(x)=c_{E}\}\\
&=\{x\in \mathcal{P}^{+}\cup\mathcal{P}^{0}\mid \Phi(x)=c_{E}\}\cup\{x\in \mathcal{P}^{-}\mid \Phi(x)=c_{E}\}\\
&\equiv L_{1}(c_{E})\cup L_{2}(c_{E}).
\end{aligned}$$
The topology of $L_{1}(c_{E})$ is simple. Indeed, by arguments before Claim 1,
it is just a manifold of codimension $1$ star shaped about the origin in $\mathcal{P}^{+}\cup\mathcal{P}^{0}$.
On the other hand, $L_{2}(c_{E})$ is a double sheet until they are glued by the nonempty set
$$\mathcal{N}_{0}=\{u\in\mathcal{P}^{-}\mid\Phi(u)=\max_{s\geq0}\Phi(su)=c_{E}\}.$$
All these topological information implies $L(c_{E})$ is not a manifold.
Notice that $\Phi$ has no critical points in $\mathcal{P}^{+}\cup(\mathcal{P}^{0}\backslash\{0\})$ and thus Claim 3 is proved.
 \hfill$\Box$

\section{Application of Theorem \ref{abstract1} to the second order case}

In this section, we focus on the periodic boundary value problem
\begin{equation}
\label{22HS}
 \ddot{x}+V'(x)=0,\quad\forall x\in\mathbb{R}^{N},\\
 \end{equation}
with
\begin{equation}
\label{periodic}
x(0)=x(T),\quad \dot{x}(0)=\dot{x}(T).
\end{equation}

In \cite{Rab1978}, Rabinowitz introduced the following variational formulation for this problem
\begin{equation}
\label{2var}
\psi(x)={1\over2}\int_{0}^{T}  |\dot{x}|^{2}dt-\int_{0}^{T}V(x)dt,\quad\forall x\in E_{T}=W^{1,2}(S_{T},\mathbb{R}^{N}),
\end{equation}
where $S_{T}=\mathbb{R}/T\mathbb{Z}$ and $W^{1,2}(S_{T},\mathbb{R}^{N})$ is endowed with the usual norm
$$||x||_{1}=\left(\int_{0}^{T}(|\dot{x}|^{2}+|x|^{2}) dt\right)^{1\over2},\quad\forall x\in E_{T}.$$
He then applied the saddle point theorem to get a nontrivial $T$-periodic solution of (\ref{22HS}).

If $V(x)=V(-x)$, there are three subspaces of $E_{T}$
$$\begin{aligned}
E_{T}^{1}=&\{x\in E_{T}\mid x\ is\ even\ about\ t = 0,\ T/2,\ and\ odd\ about\ t={T/4},\ {3T/4}\},\\
E_{T}^{2}=&\{x\in E_{T}\mid x(t+{T\over 2})=-x(t)\},\quad E_{T}^{3}=\{x\in E_{T}\mid x(-t)=-x(t)\},
\end{aligned}$$
so that the critical points of $\psi\big|_{E_{T}^{i}}$ are the ones of $\psi\big|_{E_T},$
 (cf. \cite{Long1993}, \cite{FKW2001}, \cite{WWS1996}, respectively).

Our method works for all the three subspaces, see Remark \ref{three}. Notice that $E_{T}^{1}\subset E_{T}^{2}$ and $x(t+{T\over 4})\in E_{T}^{3}$ for every $x\in E_{T}^{1}$.
We work in $E_{T}^{1}$ like  \cite{Long1993},  since in this setting a better solution of (\ref{22HS}) can be found.
First we borrow some notations and results from \cite{Long1993} and denote $E_{T}^{1}$ by $SE_{T}$.
\begin{lemma}
\label{lo93}
{\rm(Propositions 2.3 and 2.4 of \cite{Long1993})} Suppose that $V\in C^{2}(\mathbb{R}^{N},\mathbb{R}),$ then the following conclusions hold.

$1^{\circ}$ $\psi\in C^{2}(E_{T},\mathbb{R})$, i.e., $\psi$ is continuously $2$-times Fr$\acute{e}$chet differentiable on $E_{T}$.

$2^{\circ}$ There holds
$$\psi'(x)y=\int_{0}^{T}(\dot{x}\cdot\dot{y}-V'(x)\cdot y)dt,\quad\forall x,y\in E_{T}.$$

$3^{\circ}$ There holds
$$\psi''(x)(y,z)=\int_{0}^{T}(\dot{y}\cdot\dot{z}-V''(x)y\cdot z)dt,\quad\forall x,y\in E_{T}.$$

$4^{\circ}$ All the above conclusions holds if we substitute $E_T$ by $SE_T$.

$5^{\circ}$ If in addition $V(x)=V(-x)$ and $x\in SE_{T}$ is a critical point of $\psi$ on $SE_{T}$, then it is a $C^{3}$-solution of (\ref{22HS}).

\end{lemma}

By Lemma \ref{lo93}, it is enough to find a critical point of $\psi\big|_{SE_T}$ with $T$ its minimal period.
Notice that $\int_{0}^{T}xdt=0$ for every $x\in SE_{T}$. Wirtinger's inequality (\ref{20230217b}) implies
that  the usual norm $||\cdot||_{1}$ is equivalent to
$$||x||=\left(\int_{0}^{T}|\dot{x}|^{2}dt\right)^{1\over2},\quad\forall x\in SE_{T}.$$

\begin{definition}
\label{exactcpt1}
Let $S=\{e\in {SE}_{T}\mid ||e||=1\}$ be the unit sphere of $ {SE}_{T}.$ Suppose that for every $e\in {S}$, $\max_{x\in \overrightarrow{oe}} \psi(x)$ is attained and
$$c_{T}=\inf_{e\in {S}}\max_{x\in\overrightarrow{oe}} \psi(x)>0.$$
If there is a critical point $\bar{x}$ of $\psi$ on $SE_{T}$ such that
$$\psi(\bar{x})=\max_{x\in \overrightarrow{o\bar{x}}}\psi(x)=c_{T},$$
we then call $\bar{x}$ an exactly inf-max critical point of $\psi$.
 \end{definition}

Our proof of Theorem \ref{mainresult1} is based on the following simple discriminant criterion for the minimal period of an exactly inf-max critical point,
which was first observed in \cite{Xiao2010}.
\begin{lemma}
\label{criterion1}
Assume that $\bar{x}$ is an exactly inf-max critical point of $\psi$ on $SE_{T}$, then the minimal period of $\bar{x}$ is $T$.
\end{lemma}
\Proof Assume by contradiction that the minimal period of $\bar{x}$ is $T\over k$ with $k\geq 2$.
Then it can be checked directly that $x_{1}(t)=\bar{x}({t\over k})\in SE_{T}$ (cf. \cite{Xiao2010}).
 Due to the assumptions in Definition \ref{exactcpt1}, there is a constant  $\lambda_{1}>0$  such that
 $$\psi(\lambda_{1}x_{1})=\max_{x\in \overrightarrow{ox_1}}\psi(x),$$
 and by the definition of exactly inf-max critical point,
 \begin{equation}
 \label{20220728a}
 \psi(\bar{x})\leq \psi(\lambda_{1}x_{1}).
 \end{equation}

On the other hand, let $y(t)=\lambda_{1}\bar{x}(t)\in \overrightarrow{o\bar{x}}$ and then $y({t/ k})=\lambda_{1}x_{1}$. Direct computations yield that
\begin{equation}
 \label{20220728b}
 \begin{aligned}
 \psi(\lambda_{1}x_{1})=&\psi\big(y({t/ k})\big)\\
 =&\int_{0}^{T}{1\over 2k^{2}}|y'({t/ k})|^{2}dt-\int_{0}^{T}V(y({t/ k}))dt\\
 =&\int_{0}^{T/k}{1\over 2k}|y'(s)|^{2}ds-k\int_{0}^{T/k}V(y(s))ds\\
  =&\int_{0}^{T}{1\over 2k^{2}}|y'(s)|^{2}ds-\int_{0}^{T}V(y(s))ds\\
  <&\int_{0}^{T}{1\over 2}|y'(s)|^{2}ds-\int_{0}^{T}V(y(s))ds\\
 =&\psi(y).
 \end{aligned}
 \end{equation}
 By (\ref{20220728a}) and (\ref{20220728b}), we get  $$\psi(\bar{x})<\psi(y).$$ However, since $\bar{x}$ is an exactly inf-max critical point of $\psi$, we have
 $$\psi(\bar{x})=\max_{x\in \overrightarrow{o\bar{x}}}\psi(x)\geq\psi(y).$$
This leads to a contradiction and the proof is completed.\hfill$\Box$

{\bf  Proof of Theorem \ref{mainresult1} }
Due to Lemma \ref{criterion1}, we need only  find an exactly inf-max critical point of $\psi\big|_{SE_{T}}$
and so it is enough to check the conditions of Theorem \ref{abstract1} for $\psi.$

Write
$$\psi(x)={1\over2}||x||^{2}-I(x),\quad\forall x\in SE_{T},$$
where $I(x)=\int_{0}^{T}V(x)dt.$

{\bf Claim 1:} $I(x)$ is weakly continuous.

  Let $x_{n}\rightharpoonup x,$ then $x_{n}$ converges uniformly to $x$ by Proposition 1.2 of \cite{MW1989}.
  Since
$$\begin{aligned}
 |\int_{0}^{T}V(x)dt-\int_{0}^{T}V(x_{n})dt|\leq
 \int_{0}^{T}|V'(x+\theta (x -x_{n}))(x -x_{n})|dt
\end{aligned}
$$
where $0<\theta<1$, it follows immediately that $$|\int_{0}^{T}V(x)dt-\int_{0}^{T}V(x_{n})dt|\to0.$$
Thus, $\int_{0}^{T}V(x_{n})dt\to \int_{0}^{T}V(x)dt,$ and so $\int_{0}^{T}V(x)dt$ is weakly continuous.

 {\bf Claim 2:} $\lim_{x\to0}{I(x)\over||x||^{2}}=0$.

 For every  $\epsilon>0$, by (V1)  there is  $\delta>0$ such that
$$V(x)< \epsilon |x|^{2},\quad\forall |x|<\delta.$$
Since $\int_{0}^{T}xdt=0$ for $x\in SE_T$, by Sobolev inequality (see Proposition 1.3 of \cite{MW1989})
$$||x||_{\infty}^{2}\leq{T\over12}\int_{0}^{T}|\dot{x}|^{2}dt\leq {T\over12}||x||^{2}.$$
Thus for $||x||$ sufficiently small, we have $||x||_{\infty}<\delta$ and so
\begin{equation}
\label{20230217a}
V(x(t))<  \epsilon |x(t)|^{2},\quad\forall t\in[0,T].
\end{equation}

 On the other hand, by Wirtinger's inequality (see also Proposition 1.3 of \cite{MW1989})
\begin{equation}
\label{20230217b}
\int_{0}^{T}|x|^{2}dt\leq{T^{2}\over 4\pi^2}\int_{0}^{T}|\dot{x}|^{2}dt.
\end{equation}
Then by (\ref{20230217a}) and (\ref{20230217b}), we get
 $$\begin{aligned}
  \int_{0}^{T}V(x)dt\leq& \int_{0}^{T}\epsilon|x|^{2}dt
  \leq  \epsilon{T^{2}\over 4\pi^2} \int_{0}^{T}|\dot{x}|^{2}dt=\epsilon{T^{2}\over 4\pi^2}||x||^{2},
\end{aligned}$$
  and Claim 2 is proved.

  {\bf Claim 3:} ${I(sx)\over s^{2}}\to+\infty$, uniformly for $x$ on weakly compact subsets of $E\backslash\{0\}$ as $s\to+\infty$.

 As usual, we identify the equivalence class $x$ and its
continuous representant\ $$\int_{0}^{t}\dot{x}(s)ds+x(0).$$ Let $x_{n}\rightharpoonup x\neq0$, namely $x_{n}$ weakly converges to $x$.
By Proposition 1.2 of \cite{MW1989}, $x_{n}$ converges uniformly to $x$. Let $\sigma_{x}\subset[0,T]$ with positive measure $m(\sigma_{x})$ such that
 $$|x(t)|\geq 2a_{x}>0,\quad\forall t\in\sigma_{x}.$$
Thus, there is a weak open neighbourhood $O_{x}$ of $x$  such that
 $$|y(t)|\geq  a_{x}>0,\quad\forall y\in O_{x},\ t\in\sigma_{x}.$$
 Due to (V2), for every $M>0$ there is $s_{x}>0$ such that
$${V(s y(t))\over s^{2}|y(t)|^{2}}>{M\over m(\sigma_{x})a_{x}^{2}},\quad\forall s\geq s_{x},\ y\in O_{x},\ t\in\sigma_{x}.$$
As a result,
$${I(sy)\over s^{2}}=\int_{0}^{T}{V(sy)\over s^{2}}dt\geq \int_{\sigma_x} {V(sy(t))\over s^{2}}dt>M,\quad\forall s\geq s_{x},\ y\in O_{x}.$$
The rest proof is standard and we omit it.

{\bf Claim 4:} $s\mapsto{I'(sx)x\over s}$ is non-decreasing for all $x\neq0$ and $s>0.$

 By straightforward computations, we have  $ {I'(sx)x\over s}={\int_{0}^{T}V'(sx)xdt\over s},$
 and so
 $${d\over ds}\left({I'(sx)x\over s}\right)={1\over s^3}\int_{0}^{T}\Big(V''(sx)sx\cdot sx-V'(sx)sx\Big)dt\geq0,$$
 where (V3) is used in the inequality.

 As a result, there is at least one exactly $\inf$-$\max$ critical point by Theorem \ref{abstract1} and  the proof of Theorem \ref{mainresult1} is completed.\hfill$\Box$
  \begin{remark}
  \label{three}
  As we have mentioned previously that our method works for all the three subspaces $E_{T}^{1},$ $E_{T}^{2}$ and $E_{T}^3$. This is because
  $\int_{0}^{T}x(t)dt=0$ for all $x\in E_{T}^{i}$ with $i=1,2,3$. Outwardly, we can find three solutions $z_{1},z_{2}$ and $z_{3}$ with minimal period $T$.
  However, $z_1$ may be the same as $z_2$ due to $E_{T}^{1}\subset E_{T}^{2}$. In addition, it may also hold $z_{3}=z_{1}(t+{T\over4})$.
  That is, $z_{3}$ may be  the same one as $z_1$ geometrically.
  \end{remark}

\section{ Application of Theorem \ref{abstract2} to the first order case}
In this section, we use Theorem \ref{abstract2} to prove Theorem \ref{mainresult2}. For this purpose,
we first provide some preliminaries from convex analysis (cf. for instance, Chapter 2 of  \cite{Eke1990} , or Chapter 2 of \cite{MW1989}).

Let $H\in C^{1}(\mathbb{R}^{2n},\mathbb{R})$ be a convex function. Then,
$$G(y)=\sup_{x\in\mathbb{R}^{2n}}\{(x,y)-H(x)\},\quad\forall y\in\mathbb{R}^{2n},$$
is well defined and we call $G$  the Fenchel transform of $H$. It follows from the definition that $G$ is convex and lower semi continuous.

We now come to transform the conditions on $H$ in Theorem  \ref{mainresult2} to those on $G$.
\begin{lemma}
\label{convex1}
Suppose  $H\in C^{1}(\mathbb{R}^{2n},\mathbb{R})$ is strictly convex and satisfies conditions {\rm(H1)} and {\rm(H2)}. Then
$G\in C^{1}(\mathbb{R}^{2n},\mathbb{R})$ is strictly convex, $G(0)=0$ and satisfies

{\rm(G1)} $\lim_{y\to0}{G(y)\over|y|^{2}}=+\infty$,

{\rm(G2)}  There are  positive constants $R$, $b_{1}$ and $b_{2}$ such that
\begin{equation}
\label{crucialcondition2}
b_{1}|y|^{\alpha}\leq G(y)\leq b_{2}|y|^{\alpha},\quad\forall|y|\geq R,
\end{equation}
where ${1\over\alpha}+{1\over\beta}=1.$
\end{lemma}
\Proof Since  $H\in C^{1}(\mathbb{R}^{2n},\mathbb{R})$ is strictly convex, it follows from (H2) and Proposition 2.4 of \cite{MW1989} that $G\in C^{1}(\mathbb{R}^{2n},\mathbb{R})$.
Moreover, by Theorem 11.13 of \cite{Rock2009}, $G$ is strictly convex. That $G(0)=0$ follows directly from the definition.

To prove (G1), choose $\epsilon>0$ arbitrarily small. By (H1),  there is $\delta>0$ such that
 $$H(x)< \epsilon|x|^{2},\quad\forall |x|<\delta.$$
Due to the strict convexity of $H$, it holds that $$\eta=\min\{|H'(x)|,\ |x|\geq\delta\}>0.$$
Thus for every $\bar{y}\in\{y\in\mathbb{R}^{2n},\ |y|<\eta\}$, it holds  that   $|\bar{x}|=|H^{\prime-1}(\bar{y})|<\delta,$
and so
$$\begin{aligned}
G(\bar{y})=  \{(x,\bar{y})-H(x)\}\big|_{x=H^{\prime -1}(\bar{y})}
= \sup_{\{|x|\leq \delta\}}\{(x,\bar{y})-H(x)\}
\geq \sup_{\{|x|\leq \delta\}}\{(x,\bar{y})-\epsilon|x|^{2}\}
= {1\over 4\epsilon}|\bar{y}|^{2}.
\end{aligned}$$
As a result, (G1) is proved.

We now come to prove (G2). By (H2), there is $M>0$ such that
$$H(x)\leq a_{2}|x|^{\beta}+M,\quad\forall x\in\mathbb{R}^{2n}.$$
Straightforward computations yield that
\begin{equation*}
\label{20230713b}
\begin{aligned}
G(y)= \sup_{x\in\mathbb{R}^{2n}}\{(x,y)-H(x)\}
 \geq\sup_{x\in\mathbb{R}^{2n}}\{(x,y)-a_{2}|x|^{\beta}-M\}
 ={1\over\alpha}(a_{2}\beta)^{1-\alpha}|y|^{\alpha}-M.
\end{aligned}
\end{equation*}
Take $0<b_{1}<{1\over\alpha}(a_{2}\beta)^{1-\alpha}$. Then for $|y|$ large enough, we have
$G(y)\geq b_{1}|y|^{\alpha}$ and the left inequality of (\ref{crucialcondition2}) holds.
By similar arguments, one can prove the right inequality of (\ref{crucialcondition2}). \hfill$\Box$

Our proof of Theorem \ref{mainresult2} is based on the well-known Clarke dual action principle
which can be found in many literatures, for instance,  p. 411 of \cite{AM1981},
Proposition 1 of \cite{EH1985}, Theorem 2.2 of \cite{Eke1990}, and Theorem 2.3 of \cite{MW1989}. Let
$$L_{0}^{\alpha}=\left\{u\in L^{\alpha}([0,T];\mathbb{R}^{2n})\mid\int_{0}^{T}udt=0\right\}$$
and denote by $||\cdot||_{\alpha}$ the usual norm in  $L^{\alpha}([0,T];\mathbb{R}^{2n})$. Consider the functional
\begin{equation}
\label{20230811a}
\Phi(u)={1\over2}\int_{0}^{T}(Ju,\Pi u)dt+\int_{0}^{T}G(u)dt,\quad\forall u\in L_{0}^{\alpha},
\end{equation}
with
$${d\over dt}\Pi u=u\ {\rm and}\ \int_{0}^{T}\Pi u dt=0.$$
Then, if $u$ is a critical point of $\Phi$ in $L_{0}^{\alpha}$, there exists a constant\ $\xi\in\mathbb{R}^{2n}$ such that $x=J\Pi u+\xi$ is $T$-periodic solution of (\ref{1stHS}).
Therefore, we need only to find a critical point of $\Phi$ with $T$ its minimal period, in order to prove Theorem \ref{mainresult2}.

{ \bf Proof of Theorem \ref{mainresult2}.} We need only to check that $\Phi$ defined in (\ref{20230811a}) satisfies the conditions of Theorem \ref{abstract2},
with $E=L_{0}^{\alpha},$
$$a(u,u)=\int_{0}^{T}(Ju,\Pi u)dt \quad{\rm and}\quad b(u)=\int_{0}^{T}G(u)dt.$$

It is clear that $a$ is a continuous, symmetric, bilinear form. It is also well known that its positive and negative eigenspaces  are both infinite dimensional.

 Due to (G2) and the convexity of $G$, there are positive constants $c_{1}, c_{2}$ such that $$|G'(y)|\leq c_{1}|y|^{\alpha-1}+c_{2},\quad\forall y\in\mathbb{R}^{2n}.$$
Then by Lemma 4.2.8 of \cite{Eke1990}, $\Phi$ is a $C^{1}$ functional on $L_{0}^{\alpha}$.

{\bf Claim 1:} $a(u,u)$ is weakly  continuous  and $b(u)$ is weakly lower semi-continuous.

We give a proof without using the growth condition $(G2)$. Let $u_{k}$ converges weakly to $u$ in $L_{0}^{\alpha}$.
Then $u_k$ is bounded. Moreover, since $\Pi:L_{0}^{\alpha}\to L_{0}^{\beta}$ is compact (see Lemma 4.2.7 of \cite{Eke1990}),
 $u_{k}$ converges strongly to $u$ in $L_{0}^{\beta}$.
Notice that
$$\begin{aligned}
&\left|\int_{0}^{T}(Ju_{k},\Pi u_{k})dt-\int_{0}^{T}(Ju,\Pi u)dt\right|\\
\leq&\left|\int_{0}^{T}(Ju_{k},\Pi u_{k})dt-\int_{0}^{T}(Ju_{k},\Pi u)dt\right|+\left|\int_{0}^{T}(Ju_{k},\Pi u)dt-\int_{0}^{T}(Ju,\Pi u)dt\right|\\
\leq&||u_{k}||_{\alpha}||\Pi u_{k}-\Pi u||_{\beta}+\left|\int_{0}^{T}(Ju_{k},\Pi u)dt-\int_{0}^{T}(Ju,\Pi u)dt\right|.
\end{aligned}$$
Thus, $\int_{0}^{T}(Ju_{k},\Pi u_{k})dt\to\int_{0}^{T}(Ju,\Pi u)dt$ and so $a(u,u)$ is weakly continuous.

Notice that $G$ is lower semi continuous. According to Proposition 2.3.1 in \cite{Eke1990}, $b(u)$ is lower semi continuous  on $L_{0}^{\alpha}$.
On the other hand, since $G(u)$ is convex, $b(u)$ is also convex.
By Theorem 1.2 in \cite{MW1989}, it is weakly lower semi-continuous.
As a result, $b(u)$ is weakly lower semi-continuous on $L_{0}^{\alpha}$.

{\bf Claim 2:} $b(0)=0$, $\lim_{u\to0}{b(u)\over||u||_{\alpha}^{2}}=+\infty $ and $\lim_{u\to\infty}{b(u)\over||u||_{\alpha}^{2}}=0.$

 The facts that $b(0)=0$ and $\lim_{u\to\infty}{b(u)\over||u||_{\alpha}^{2}}=0$  follow  immediately from $G(0)=0$ and (G2), respectively.
 Now we come to prove
 \begin{equation}
 \label{20230811b}
 \lim_{u\to0}{b(u)\over||u||_{\alpha}^{2}}=+\infty.
 \end{equation}

For $M=2M_{1}T^{2-\alpha\over\alpha}$ with $M_{1}>0$ arbitrarily large, condition (G1) implies that there is sufficiently small $\delta>0$ such that
\begin{equation}
\label{83b}
G(y)\geq M|y|^{2},\quad\forall |y|\leq\delta.
\end{equation}
 Consider the set $\{y\in\mathbb{R}^{2n}\mid |y|\geq\delta\}.$
 On the one hand, by (G2) we have $$G(y)\geq b_{1}|y|^{\alpha},\quad \forall |y|\geq R.$$
 On the other hand, the fact that $G(y)>0$ for $y\neq0$ implies the existence of $d>0$ small enough so that
 $$G(y)\geq d|y|^{\alpha},\quad\forall \delta\leq|y|\leq R.$$
 In summary, we get
 \begin{equation}
\label{83c}
G(y)\geq d |y|^{\alpha},\quad\forall |y|\geq\delta.
\end{equation}

Write $u=v+w$ with
$$v(t)=\left\{\begin{array}{ll}
u(t),&\quad{\rm if}\ |u(t)|\leq\delta,\\
0,&\quad{\rm if}\ |u(t)|> \delta.\\
\end{array}
\right.$$

On the other hand, by (\ref{83b}), (\ref{83c}) and H$\ddot{\rm o}$lder inequality it follows that
\begin{equation}
\label{83e}
\begin{aligned}
 \int_{0}^{T}G(u)dt&=\int_{0}^{T}G(v)dt+\int_{0}^{T}G(w)dt\\
 &\geq M\int_{0}^{T}|v|^{2}dt+d\int_{0}^{T}|w|^{\alpha}dt\\
 &=M||v||_{2}^{2}+d||w||_{\alpha}^{\alpha}\\
 &\geq MT^{-{2-\alpha\over\alpha}}||v||_{\alpha}^{2}+d||w||_{\alpha}^{\alpha}\\
 &= 2M_{1}||v||_{\alpha}^{2}+d||w||_{\alpha}^{\alpha}.
 \end{aligned}
\end{equation}

Notice $1<\alpha<2$ and $||w||_{\alpha}\leq||u||_{\alpha}.$ If $||u||_{\alpha}$ is sufficiently small, it holds
$$d||w||_{\alpha}^{\alpha} \geq 2M_{1}||w||_{\alpha}^{2}.$$
Therefore, for $||u||_{\alpha}$ small enough we have
$$b(u)\geq 2M_{1}(||v||_{\alpha}^{2}+||w||_{\alpha}^{2})\geq M_{1}( ||v||_{\alpha}+||w||_{\alpha})^{2}\geq M_{1}||v+w||_{\alpha}^{2}=M_{1}||u||_{\alpha}^{2}.$$
That is, (\ref{20230811b}) holds.

{\bf Claim 3:} $b'(u)u\geq b(u)>0$ for $u\neq0$.

Since $G(0)=0$ and $G$ is strictly convex, for $y\neq0$ it holds that $G(y)>0$  and $G'(y)y>G(y)$. Then the former inequality follows immediately.

{\bf Claim 4:}  For every $u\neq 0$, $s\mapsto{b'(su)u\over s}$ is  nonincreasing on $(0,+\infty)$.

Since ${b'(su)u\over s}=\int_{0}^{T}{G'(su)u\over s}dt$, Claim 4 follows from (H3).

All the conditions of Theorem \ref{abstract2} are checked and so there is a critical point $\bar{u}\in\mathcal{P}^{-}$ such that
$$\Phi(\bar{u})=\max_{s\geq0}\Phi(s\bar{u})=\inf_{e\in\mathcal{S}^{-}}\max_{s\geq0}\Phi(se).$$
The fact that the minimal period of $\bar{u}$ is $T$ follows by the same arguments as those in Lemma \ref{criterion1}. The only thing worth explaining here is that
$\bar{u}({t\over k})\in\mathcal{P}^{-}.$ This is because $$a\left(\bar{u}\left({t\over k}\right),\bar{u}\left({t\over k}\right)\right)=ka(\bar{u}(t),\bar{u}(t))<0.$$
\hfill$\Box$

Finally, we come to explain how to use our idea to study the Nehari-Pankov manifold with its application to the minimal period problem for general Hamiltonian systems.
For convenience of writing, we focus on the variational setting of the second order Hamiltonian systems without the even condition on $V$.

Let $$SE_{T}=\left\{x\in W^{1,2}([0,T],\mathbb{R}^{N})\mid x(t)=x(-t)\right\}$$ with the usual norm $||\cdot||$.
 Write
 $SE_{T}=\mathbb{R}^{N}\oplus\widetilde{SE_{T}}$ with
  $\widetilde{SE_{T}}=\{x\in SE_{T}\mid\int_{0}^{T}xdt=0\}.$
It is well known that the critical point of the functional
$$\psi(x)={1\over2}\int_{0}^{T}|\dot{x}|^{2}dt-\int_{0}^{T}V(x)dt,\quad\forall x\in SE_{T},$$
is a $T$-periodic solution of (\ref{2ndHS}).

If for every $e\in\widetilde{S}= \widetilde{SE_{T}}\cap\{x\in SE_{T}\mid||x||=1\}$,
there is a unique maximum point $\mathcal{M}_{e}$ of $\psi$ on the half-space $H_{e}=\{\lambda e\mid\lambda\geq0\}\oplus\mathbb{R}^{N},$ we can define the Nehari-Pankov manifold
$$\mathcal{M}=\{\mathcal{M}_{e}\mid e\in\widetilde{SE_{T}}\}.$$
Indeed, our subsequent arguments are valid at least if every $\mathcal{M}_{e}$ is contractible.

Define
$$c_{T}=\inf_{x\in\mathcal{M}}\psi(x)=\inf_{e\in\widetilde{S}}\max_{x\in H_{e}}\psi(x).$$
It can be proved that $c_{T}>0$ can be achieved under the usual superquadratic conditions. Like in section 2, define
$\mathcal{M}_{+}=\{x\in\mathcal{M}\mid \psi(x)>c_{T}\}$ and $\mathcal{M}_{0}=\{x\in\mathcal{M}\mid \psi(x)=c_{T}\}$.
Consider the level set $L(c_{T})$. Then for every $x\in \mathcal{M}_{+}$, $L(c_{T})\cap H_{x}$ is an $N$-dimensional sphere.
While if $x\in\mathcal{M}_{0}$, $L(c_{T})\cap H_{x}=\{x\}.$ As a result, $L(c_{T})$ can not be a manifold and so there is at least one critical point $\bar{x}\in\mathcal{M}_{0}$.
Then one can follow the line of Lemma \ref{criterion1} to prove that $T$ is the minimal period of $\bar{x}.$

However, we do not know how to ensure the unique condition on maximum points of $\psi|_{H_e}$. Even the convexity assumption is not enough.
This is in sharp contrast to some functionals for elliptic equations (cf. Proposition 39 in \cite{Szulkin-Weth2010}).

\noindent{\bf Acknowledgements} The authors would like to sincerely appreciate Professor Shiqing Zhang for many helpful discussions on convex analysis.

\end{document}